\documentclass[a4paper,12pt,twoside]{article}
\usepackage{amssymb}
\usepackage[final]{lpretex}
\usepackage{title}
\usepackage{a4}
\usepackage{amscd}
\input xy
\xyoption{all}

\EndOfDump

\def\DHpartformat#1{(#1) }

\LBsetstyleof{partno}{arabic}

\DocVersion[Changed to LaTeX]{1}

\let\define\def

  \def\C {{\mathbb C}}

  \def\P {{\mathbb P}} 
\def\Q {{\mathbb Q}} 
\def\T {{\mathbb T}} 
  \def\X {{\mathbb X}}
\def\Z {{\mathbb Z}} 

\define \n {\mathbb N}
\define \z {\mathbb Z}
\define \q {\mathbb Q}
\define \PP {\mathbb P}

\def\sA {{\Cal A}}  
 \def\sE {{\Cal E}} \def\sF {{\Cal F}}

 \def\sN {{\Cal N}} \def\sO {{\Cal O}}
 
 \def\sT {{\Cal T}} \def\sU {{\Cal U}}
  \def\sX {{\Cal X}}
\def\sY {{\Cal Y}}

\define \cN {\Cal N}
\define \cf {\Cal F}
\define \cg {\Cal G}
\define \cE {\Cal E}
\define \ce {\Cal E}
\define \cc {\Cal C}
\define \cV {\Cal V}
\define \cA {\Cal A}
\define \cK {\Cal K}
\define \cO {\Cal O}
\define \cF {\Cal F}
\define \cn {\Cal N}
\define \cI {\Cal I}
\define \sP {\Cal P}
\define \sW {\Cal W}
\def\tA {\widetilde{\Cal A}}

\def\s {\sigma}

\define \x {\xi}
\define \y {\eta}
\define \G {\Gamma}
\define \r {\rho}
\define \w {\omega}

\def \tU {\widetilde U}

\def\tX {\widetilde X}

\def \trho {\widetilde {\rho}}

\def \tp {\widetilde{\mathbb P}}
\define \tH {\widetilde H}
\define \tG {\widetilde{\Gamma}}
\define \tW {\widetilde W}
\define \tF {\widetilde F}
\define \tm {\widetilde m}
\define \St {\widetilde S}
\define \Xt {\widetilde X}
\define \tS {\widetilde S}
\define \tpsi {\widetilde \psi}
\define \tL {\widetilde L}
\define \tE {\widetilde E}
\define \tl {\widetilde l}
\define \tA {\widetilde A}
\define \tom {\widetilde\omega}
\define \tT {\widetilde T}
\define \tB {\widetilde B}
\define \tf {\widetilde f}
\define \tsA {\widetilde{\sA}}

\define \tM {\widetilde M}
\define \tphi {\widetilde{\phi}}
\define \trho {\widetilde{\rho}}
\define \tR {\widetilde R}
\define \tp {\widetilde p}
\define \tq {\widetilde q}
\define \tc {\widetilde c}
\define \tsF {\widetilde {\sF}}
\define \tsN {\widetilde {\sN}}
\define \tsU {\widetilde {\sU}}
\define \tsX {\widetilde {\sX}}
\define \th {\widetilde h}

\def\pd {\partial}

\def \Dx1 {\frac{\pd}{{\pd} x_1}}
\def \Dy1 {\frac{\pd}{{\pd} y_1}}
\def \Dz1 {\frac{\pd}{{\pd} z_1}}
\def \Dx2 {\frac{\pd}{{\pd} x_2}}
\def \Dy2 {\frac{\pd}{{\pd} y_2}}
\def \Dz2 {\frac{\pd}{{\pd} z_2}}

\def\q {\quad} 

\def\mapdiagr#1{\Big\searrow\rlap{$\raise 5pt\vbox{{\hbox{$\mkern -15mu\scriptstyle#1$}}}$}}   

\def\mapdiagl#1{\llap{$\raise 5pt\vbox{{\hbox{$\scriptstyle#1\mkern
-15mu$}}}$}\Big\swarrow}              

\def\Mapdiagr#1{\nearrow\rlap{$\lower 5pt\vbox{{\hbox{$\mkern
-15mu\scriptstyle#1$}}}$}} 

\def\Mapdiagl#1{\llap{$\lower 5pt\vbox{{\hbox{$\scriptstyle#1\mkern
-15mu$}}}$}\searrow} 

\def\Mapswr#1{\swarrow\rlap{$\lower 5pt\vbox{{\hbox{$\mkern
-15mu\scriptstyle#1$}}}$}}              

\def\Mapnwl#1{\nwarrow\rlap{$\lower 5pt\vbox{{\hbox{$\mkern
-15mu\scriptstyle#1$}}}$}}

\def \inj {\hookrightarrow}

\define \Rhook {\hookrightarrow}

\def \half {\raise1pt\hbox{$\scriptstyle
        \frac{1}{2}\displaystyle$}}

\def \x{{\sl X}\llap{$\mkern -2mu {\scriptstyle -}$}}
\def \Symm {\operatorname{Sym}}

\define \Kod {\operatorname{Kod}}
\define \dimension {\operatorname{dim}}
\define \codim {\operatorname{codim}}
\define \contr {\operatorname{contr}}
\define \rk {\operatorname{rank}}
\define \im {\operatorname{im}}
\define \Mor {\operatorname{Mor}}
\define \Cl {\operatorname{Cl}}
\define \Hilb {\operatorname{Hilb}}
\define \degree {\operatorname{deg}}
\define \mult {\operatorname{mult}}
\define \Aut {\operatorname{Aut}}
\define \NS {\operatorname{NS}}
\define \Gal {\operatorname{Gal}}
\define \ch {\operatorname{char}}
\define \Jac {\operatorname{Jac}}
\define \Km {\operatorname{Km}}
\define \Sec {\operatorname{Sec}}
\define \Stab {\operatorname{Stab}}
\define \Br {\operatorname{Br}}
\define \inv {\operatorname{inv}}
\define \tr {\operatorname{tr}}
\define \Frob {\operatorname{Frob}}
\define \Symn {\operatorname{Sym}^n}
\define \Ev {\sE^\vee}
\define \ordp {\operatorname{ord}_p}
\define \Supp {\operatorname{Supp}}
\define \Ann {\operatorname{Ann}}
\define \disc {\operatorname{disc}}
\define \Lie {\operatorname{Lie}}
\define \embdim {\operatorname{embdim}}

\def\hod#1#2#3#4{\ensuremath{\if#30 H^{#2}({#1},{\cal O}_{#1}) \else 
 H^{#2}(#1,\Omega^{#3}\if\relax{#4}\relax_{#1}\else _{#1/#4}\fi)\fi}}

\begin{document}
\title[Terminal singularities]{The singularities of 
$A_{\MakeLowercase {g}}^P$}

\author{J. Armstrong}
\author{N. I. Shepherd-Barron}
\email{john.armstrong@kcl.ac.uk, nicholas.shepherd-barron@kcl.ac.uk}
\address
{Math. Dept.\\
King's College London\\
Strand\\
WC2R 2LS\\
UK}
\maketitle
\medskip
Let $\sA_g$ denote the stack, and $A_g=[\sA_g]$ the coarse
moduli space,
of principally polarized abelian $g$-folds
and $\sA_g^P$, $A_g^P=[\sA_g^P]$ their perfect cone compactifications;
these are particular toroidal compactifications.

The point of this paper is to complete the proof in [SB] of 
a result that is suggested by Tai's  
paper \cite {T}, that in characteristic zero the singularities of 
$A_g^P$ are canonical if $g\ge 5$ and terminal if
$g\ge 6$. However, Hulek and Sankaran have pointed out
that the argument given in [SB] is incomplete, because Proposition 3.2 of [SB]
is not strong enough. We complete it by replacing 
Proposition 3.2 of \emph{loc. c it.} with Theorem \ref{canonical} below.
\begin{section}{The proof}
We recall the RST criterion for the geometric
quotient $[\sX]$ of a smooth Deligne--Mumford stack $\sX$
over a field of characteristic zero 
to have canonical or terminal singularities.
We phrase it in terms of the action of a finite group
$G$ acting on a smooth variety $Y$.
First,
let $\{x\}$ denote the fractional
part of a real number $x$. Then,
for any $s\in G-\{1\}$,
define $\lambda_{Y,P}(s)=\sum\arg(\zeta)/2\pi$, where
the sum is over the eigenvalues $\zeta$ of $s$ on 
the tangent space $T_{Y,P}$.

\medskip
(RST) $G$ acts freely in codimension one in a neighbourhood of $P$ and
$[Y/G]$ has canonical (resp., terminal)
singularities at the image of $P$ if and only if 
$\lambda_{Y,P}(s)\ge 1$ (resp., $\lambda_{Y,P}(s)>1$)
for every $s\in  G-\{1\}$ and every choice of $\zeta$.
\medskip

\begin{theorem}\label{canonical} 
The singularities of $A_g^P$ are canonical if $g= 5$, and terminal if
$g\ge 6$.
\begin{proof}
The problem reduces to this. Suppose that $h,r\ge 0$, $g=h+r$,
$C$ is a principally polarized abelian $h$-fold, $W=H^0(C,\Omega^1_C)^\vee$,
$\Lambda=\Z^r$, $V=\Symm^2W\oplus (W\otimes\Lambda)$
and $\T$ is the torus with $\X_*(\T)=B(\Lambda)$, the group of
symmetric $\Z$-valued
bilinear forms on $\Lambda$. There is a torus embedding $\T\inj X$
with terminal singularities and a $\T$-bundle $\sT\to V$ such that
locally $\sA_g^P=\sX/G$ for a finite group $G=\Stab_{\sX}(P)$,
where $\sX=\sT\times^\T X\to V$ is a $G$-equivariant $X$-bundle
and $P\in\sX$. This is compatible with a $\Z$-linear action of $G$ on
$\Lambda$ and an action of $G$ on $C$, so that $\T\inj X$
is $G$-equivariant and $G$ acts linearly on $V$. In terms of
co-ordinates, choose a $\Z$-basis $(x_1,\ldots,x_r)$ of $\Lambda$
and write $Z_{pq}=\exp 2\pi i(x_p\odot x_q)$. Then locally
$\sO_{\tsX}=\sO_V[Y_m]$ for some monomials $Y_m$ in the $Z_{pq}$
and the action of an element $s\in G$ is given by
$s^*Y_m=\Phi_{s,m}W_{s(m)}$, where $\Phi_m\in\sO^{hol}_V$,
$\Phi_m$ does not vanish at $0$ and $W_{s(m)}$ is a monomial
in the $Z_{pq}$.


\begin{lemma} If $h\ge 1$, then $G$ acts effectively on $V$.
\begin{proof} The only elements of $GL(W)\times GL(\Lambda)$
that act trivially on $V$ are $\pm 1$.
\end{proof}
\end{lemma}

The next proposition slightly extends one of Tai's. We omit the proof, but it 
follows his almost exactly: reduce to the case where the order of $s$
divides $12$ and then check a finite list. We
executed this check by writing a routine in Mathematica.

\begin{proposition}\label{crux} Assume that $g\ge 5$, that $h\ge 1$
and that $s\in G-\{1\}$.

\part[i] If $s\vert_W\ne\pm 1$, then
$\lambda_{\Symm^2W}(s)\ge 1$ 
if $h=5$ and $\lambda_{\Symm^2W}(s)>1$ if $h\ge 6$.
\part[ii] If $s\vert_V$ is not of order $2$,
then $\lambda_V(s)\ge 1$. 
\part[iii] If $\lambda_V(s)<1$ then $h=1$, $s\vert_W=-1$,
$s\vert_{\Lambda_{\Q}}=(1,(-1)^{r-1})$ and $\lambda_V(s)=1/2$.
\noproof
\end{proposition}

In particular, taking $g=h$ shows that the singularities of the
interior $A_g$ are terminal if $g\ge 6$ and canonical if $g=5$.
 
For the proof of the theorem,
choose a $\T\rtimes G$-equivariant resolution $\tX\to X$, and put
$\tsX=\sT\times^{\T}\tX$. 
If $s$ acts trivially on some divisor $E$ in
$\tsX$, then $E$ is in the boundary $\partial\tsX=\tsX-\sT$. 
Each boundary divisor dominates $V$, since $\tsX\to V$
is a fibre bundle, and then $s$ acts trivially on $V$. This is false,
so that $G$ acts freely in codimension one on $\tsX$.

For some $q\ge 1$ there is a $G$-invariant 
generator $\s$ of the sheaf $\sO(qK_{\sX})$. Then
$\s$ is also a generator of $\sO(qK_{[\sX/G]})$.
Since $\sX$ is smooth, $\s$ vanishes along all exceptional divisors in $\tsX$,
so that $\s$ also vanishes along all exceptional divisors in $[\tsX/G]$.
So it is enough to deal with the singularities of $[\tsX/G]$.

Suppose that $h\ge 1$ and $g\ge 6$.
Assume that $\lambda_{\tsX,P}(s)\le 1$.
Note that $\lambda_{\tsX,P}(s)\ge\lambda_V(s)+\lambda_{\tsX/V,P}(s)$,
from the definition of $\lambda$,
and that $P$ lies in $\partial\tsX$,
since $\lambda_{\sA_g,Q}>1$ for all points $Q$ of $\sA_g$.
If $\lambda_{\tsX/V,P}(s)=0$ then 
$s$ acts trivially on the relative tangent space $T_{\tsX/V,P}$,
so that a suitable point $Q\in\sT$ 
that is near $P$ would give $Q\in\sA_g$ with
$\lambda_{\sA_g,Q}(s)=\lambda_{V,P}(s)\le 1$. 
So $\lambda_{\tsX/V,P}(s)>0$, $\lambda_V(s)<1$, and then
$s$ is of order $2$, $\lambda_V(s)=1/2$,
$h=1$ and $r\ge 5$. Moreover, $\lambda_{\tsX/V,P}(s)=1/2$.

There is an $s$-equivariant factorization
$\tsX\stackrel{\alpha}{\to}\sY\stackrel{\beta}{\to} V$
of $\tsX\to V$, where $\alpha$ is 
a vector bundle
and $\beta$ is a torus bundle. The algebra $\sO_{\sY}$ is generated 
as an $\sO_V$-algebra by
those monomials $Y_m$ that are invertible at $P$.

If $\lambda_{\tsX/\sY,P}(s)=0$, then we could,
as before, find $Q\in\sA_g$
with $\lambda_{\sA_g,Q}(s)=1$, contradiction. 
So $\lambda_{\tsX/\sY,P}(s)=1/2$ and $\lambda_{\sY/V,P}(s)=0$. 
We have $\sO_{\tsX}=\sO_{\sY}[Y_1,\ldots,Y_n]$
and $s^*Y_i=\Psi_iY_{s(i)}$, with $\Psi_i\in\sO_{\sY}^*$.
Since $\lambda_{\tsX/\sY,P}(s)=1/2$,
the trace of $s$ acting on $T_{\tsX/\sY,P}$ equals
$\dim T_{\tsX/\sY,P}-2$, so that
there at most two indices $i$ for which $s(i)\ne i$.
Since also $\lambda_{\sY/V,P}(s)=0$, there is a $\Q$-basis of
$B(\Lambda_{\Q})$ of which $s$ fixes at least ${{r+1}\choose 2}-2$ 
elements. But $s\vert_{\Lambda\otimes\Q}=(1,(-1)^{r-1})$ and $r\ge 5$,
so this is impossible and we have contradicted the assumption
that $\lambda_{\tsX,P}\le 1$.

The proof that $A_5^P$ has canonical
singularities when $h\ge 1$ is the same.

Finally, if $h=0$, then we use a theorem of Snurnikov [Sn].
Since it remains unpublished, we include a proof.


\begin{proposition} Suppose 
 $\T\inj X$ is a $G$-equivariant
torus embedding, $U\subset X$ a $G$-invariant neighbourhood
of $X-\T$, $G$ acts freely in codimension one
on $U\cap\T$ and that $U$ and
$[(U\cap\T)/G]$ have terminal (canonical) singularities. Then $G$ acts freely in
codimension one on $U$ and
$[U/G]$ has terminal (canonical) singularities.
\begin{proof} 
Let $\pi:\tX\to X$ be a $\T\rtimes G$-equivariant resolution. 
As before, it is enough to show that the $G$-action on $\tU=\pi^{-1}(U)$
satisfies (RST).

If $s\in G$ and $s(P)=P$, then, locally near $P$, $\tX$ is an
$s$-equivariant vector bundle 
$q:\tX=\Sp\C[x_,y_j^\pm]\to\T_1=\Sp\C[y_j^\pm]$;
the $y_j$ are those characters of $\T$
that are invertible at $P$ and $s^*x_i=\phi_i x_{s(i)}$,
where $\phi_i$ is a monomial in the $y_j$.
Let $R$ denote the origin in the fibre $F=q^{-1}(q(P))$;
then $\lambda_{\tU,P}(s)=\lambda_{\tU,R}(s)$,
since $s\mapsto \lambda(s)$ is constant on the line in $F$
through $R$ and $P$.

So we can assume that $P=R$ and that $\tX\ne\T_1$. 
Then the locus $x_i=t$,
$y_j^\pm =1$ for all $i,j$ is a line $L$ of $s$-fixed
points in $\tX$ that passes through $P$ and meets $\tU\cap\T=U\cap \T$.
Suppose $Q\in U\cap L$ is general; then
$\lambda_{\tX,P}(s)=\lambda_{U\cap\T,Q}(s)$.

\end{proof}
\end{proposition}
This completes the proof of the theorem.
\end{proof}
\end{theorem}

\end{section}
\bibliography{alggeom,ekedahl}
\bibliographystyle{pretex}
\end{document}